\documentclass{article}
\usepackage{amsmath, amsthm, amssymb}
\usepackage{makeidx}
\usepackage{amsbsy}
\usepackage{graphicx}
\begin{document}
% generates the title
\def\theequation{\arabic{section}.\arabic{equation}}
\setcounter{page}{1}

\title{\bf Estimation of mean vector in elliptical models}

\vspace{3.5cm}

\author{{ M. Arashi}\vspace{.5cm}
\\\it Faculty of Mathematics, Shahrood University of Technology
\\\it Shahrood, Iran\\{\em Email: m\_arashi\_stat@yahoo.com}}

\date{}
\maketitle

\begin{quotation}
\noindent {\it Abstract:} In this paper, we are basically
discussing on a class of Baranchik type shrinkage estimators of
the vector parameter in a location model, with errors belonging to
a sub-class of elliptically contoured distributions. We derive
conditions under Schwartz space in which the underlying class of
shrinkage estimators outperforms the sample mean. Sufficient
conditions on dominant class to outperform the usual James-Stein
estimator are also established. It is nicely presented that the
dominant properties of the class of estimators are robust truly
respect to departures from normality.
\par

\vspace{9pt} \noindent {\it Key Words and Phrases:} Elliptically
contoured distribution, James-Stein estimator, Jeffreys' prior,
Minimax, Robustness, Shrinkage estimator, Schwartz space.
\par
\end{quotation}\par

\newcommand{\by}{\boldsymbol{y}}
\newcommand{\bz}{\boldsymbol{z}}
\newcommand{\ba}{\boldsymbol{a}}
\newcommand{\bY}{\boldsymbol{Y}}
\newcommand{\bb}{\boldsymbol{b}}
\newcommand{\bx}{\boldsymbol{x}}
\newcommand{\bX}{\boldsymbol{X}}
\newcommand{\bA}{\boldsymbol{A}}
\newcommand{\bB}{\boldsymbol{B}}
\newcommand{\bC}{\boldsymbol{C}}
\newcommand{\bH}{\boldsymbol{H}}
\newcommand{\bM}{\boldsymbol{M}}
\newcommand{\bR}{\boldsymbol{R}}
\newcommand{\bS}{\boldsymbol{S}}
\newcommand{\bI}{\boldsymbol{I}}
\newcommand{\bSi}{\boldsymbol{\Sigma}}
\newcommand{\0}{\boldsymbol{0}}
\newcommand{\bth}{\boldsymbol{\theta}}
\newcommand{\bet}{\boldsymbol{\eta}}
\newcommand{\be}{\boldsymbol{\epsilon}}
\newcommand{\bPh}{\boldsymbol{\Phi}}
\newcommand{\bd}{\boldsymbol{\delta}}
\newcommand{\W}{\mathcal{W}}
\newcommand{\A}{\mathcal{A}}
\newcommand{\B}{\mathcal{B}}
\newcommand{\E}{\mathcal{E}}
\newcommand{\N}{\mathcal{N}}
\newcommand{\D}{\mathcal{\boldsymbol{D}}}
\newcommand{\beq}{\begin{eqnarray}}
\newcommand{\eeq}{\end{eqnarray}}
\newtheorem{thm}{Theorem}[section]
\newtheorem{cor}{Corollary}[section]
\newtheorem{lem}{Lemma}[section]
\newtheorem{rem}{Remark}[section]

\section{Introduction}
The assumption of normality restricts the range of possible applications, especially in flatter densities. The elliptically contoured distributions (ECDs) are the parametric forms of the spherical symmetric distributions, which are invariant under orthogonal transformations and have equal density on sphere if densities exist. ECD's primary purpose is to provide a highly impressive list of heavier/lighter tail alternatives to the multivariate Gaussian models. Materials involving vector-variate distributional properties and inferential problems will be found entirely in the work of a couple of statisticians like Das Gupta et al. (1972), Cambanis et al. (1981), Muirhead (1982), Anderson et al. (1986), Cellier et al. (1989), Anderson and Fang (1990), Fang and Zhang (1990), Fang et al. (1990) and Kibria and Haq (1999). Among others, the book of Gupta and Varga (1993) illustrates some significant results dealing with matrix-variate ECD. In addition, some of the well-known elliptical distributions are the Gaussian, Pearson Type II/VII,
 Student's t, logistics, Kotz type, Laplace, Bessel and power exponential multivariate distributions.

In this paper, we consider the location model in a more general
setup involving dependent errors. Initially let $\mathcal{S}(p)$
denotes the set of all $p\times p$ positive definite matrices. The
precise set-up of the problem is as follows: Let $\bY_i$ be an
$p\times1$ response vector with model \beq\label{eq11}
\bY_i&=&\bth+\be_i, \ \ \ 1\leq i\leq N. \eeq Here $\bth$ is a
$p\times1$ vector of location parameters and $\be_i$ is a
$p\times1$ error vector such that \beq\label{eq12} E(\be_i)=\0, \
\ Cov(\be_i\be_j)=\bSi\in\mathcal{S}(p), \ \ \ i,j=1,\cdots,N, \
\ \ N>p. \eeq
 It is assumed, in general, $\be=(\be_1,\cdots,\be_N)'$ have a joint elliptically contoured distribution. Typically if it possess a density, it is followed by
\beq\label{eq13}
f(\be|\bSi)\propto|\bSi|^{\frac{-N}{2}}g\left(tr\;\bSi\sum_{i=1}^N
\be_i\be_i^\prime\right), \eeq where $g(.)$ is a non-negative
function over $\mathbb{R}_+$ such that $f(.|.)$ is a density
function with respect to (w.r.t.) a $\sigma$-finite measure $\mu$ on $\mathbb{R}^p$. In this case, notation $\be_i\sim\E_{p}(\0,\bSi,g)$ would probably be used.

Due to Chu (1973), each component of the aforementioned model being proposed in
\eqref{eq13}, possibly can be presented as the following form.
\beq\label{eq14}
f_{\be_i}(\bx)=\int_0^\infty
\W(t)\phi_{\N_p(\0,t^{-1}\bSi)}(\bx)dt, \eeq where
$\phi_{\N_p(\0,t^{-1}\bSi)}(.)$ is the $pdf$ of
$\N_p(\0,t^{-1}\bSi)$, \beq\label{eq15}
\W(t)=(2\pi)^{\frac{p}{2}}|\bSi|^{\frac{1}{2}}t^{-\frac{p}{2}}\mathcal{L}^{-1}[f(s)],
\eeq $\mathcal{L}^{-1}[f(s)]$ denotes the inverse Laplace
transform of $f(s)$ with $s=t[\bx'\bSi^{-1}\bx/2]$.

The inverse Laplace transform of $f(.)$ exists provided that the following conditions are satisfied.
\begin{enumerate}
\item[(i)]$f(t)$ is differentiable when $t$ is sufficiently large.
\item[(ii)]$f(t)=o(t^{-m})$ as $t\rightarrow\infty$, $m>1$.
\end{enumerate}
Although, it is rather difficult to derive the inverse Laplace transform of some functions, we are able to handle it for many density generators of elliptical densities. Refer to Debnath and Bhatta (2007) for more specific details.

The mean of $\be_i$ is the zero-vector and the covariance-matrix of
$\be_i$ is
\beq\label{eq16}
Cov(\be_i)&=&\int_0^\infty Cov(\be_i|t)\W(t)dt\nonumber\\
&=&\int_0^\infty \W(t)Cov\left\{\N_p(\0,t^{-1}\bSi)\right\}dt\nonumber\\
&=&\left(\int_0^\infty t^{-1}\W(t)dt\right)\bSi,
\eeq
provided the above integral exists.

Another sub-class of ECDs which includes the above class may be
generated by a signed measure $\W$ on the measurable space
$(\mathbb{R}^+,\mathbb{B})$ such that the pdf $f(.)$ can be
expressed through the following procedures:
\beq\label{eq17}
(i)&&f(\bx)=\int_0^\infty \phi_{\N_p(\0,t^{-1}\bSi)}(\bx) \W(dt),\\
(ii)&&\int_0^\infty t^{-1}\W^+(dt)<\infty,\nonumber\\
(iii)&&\int_0^\infty t^{-1}\W^-(dt)<\infty,\nonumber \eeq where
$\W^+ - \W^-$ is the Jordan decomposition of $\W$ in positive and
negative parts (see e.g. Srivastava and Bilodeau, 1989).
Note that from $(ii)-(iii)$ of \eqref{eq17},
\beq\label{eq18}
\int_0^\infty t^{-1}\W(dt)<\infty
\eeq
and thus, $Cov(\be_i)$
exists under the sub-class defined above.

Now, under Bayesian framework, it is properly assumed that in distribution, little is known ofcourse, a priori, about
the parameters, the elements of $\bth$ and the $p(p+1)/2$
distinct elements of $\bSi\in\mathcal{S}(p)$. We shall first of
all suppose that the elements of $\bth$ and those of $\bSi$ are
approximately independent (see Box and Tiao, 1992, page 425), i.e.
\begin{eqnarray}\label{eq19}
\pi(\bth,\bSi)\doteq\pi(\bth)\pi(\bSi).
\end{eqnarray}
Using the invariant theory due to Jeffreys (1961), we take
\begin{eqnarray*}
\pi(\bth)&\propto&\mbox{constant},\nonumber\\
\pi(\bSi)&\propto&|\bSi|^{-\frac{p+1}{2}},
\end{eqnarray*}
as the prior knowledge about the parameter space.\\
Next step being taken, is giving results for the marginal posterior distribution of the location parameter given responses.
\begin{lem}
Assume in the location model \eqref{eq11}, $\be_i\sim
\E_p(\0,\bSi,g)$, where $\bSi\in\mathcal{S}(p)$. Then, w.r.t. the prior distribution given by \eqref{eq19}, the
posterior distribution of $\bth$ is multivariate Student's t
distribution, denoted by $\bth|\bY\sim t_p(\bar{\bY},\bS,N-p)$,
with the following $pdf$
\begin{eqnarray*}
f(\bth|\bY)=\frac{N^{\frac{p}{2}}\Gamma\left(\frac{N}{2}\right)|\bS|^{-\frac{1}{2}}}
{\pi^{\frac{n}{2}}(N-p)^{\frac{N-p}{2}}\Gamma\left(\frac{N-p}{2}\right)}\bigg[1+N(\bth-\bar{\bY})'\bS^{-1}(\bth-\bar{\bY})\bigg]^{-\frac{N}{2}}
\end{eqnarray*}
where $\bY=(\bY_1,\cdots,\bY_N)$, and \beq\label{eq110}
\bar{\bY}=\frac{1}{N}\sum_{i=1}^N\bY_i,\hspace{.5cm}\bS=\sum_{i=1}^N(\bY_i-\bar{\bY})(\bY_i-\bar{\bY})'.
\eeq
\end{lem}
\noindent\textbf{Proof:} Using Proposition 1 of Ng (2002), one
can directly obtain \beq\label{eq111}
f(\bth|\bY)\propto\left|\sum_{i=1}^N(\bY_i-\bth)(\bY_i-\bth)'\right|^{-\frac{N}{2}},
\eeq
which is the same as we take the errors to be normally distributed (Zellner, 1971, P.243).\\
At this level, through making conclusion based on the following equality
\begin{eqnarray*}
(\bY_i-\bth)(\bY_i-\bth)'&=&(\bY_i-\bar{\bY})(\bY_i-\bar{\bY})'+(\bar{\bY}-\bth)(\bar{\bY}-\bth)'\\
&&+2(\bY_i-\bar{\bY})(\bar{\bY}-\bth)',
\end{eqnarray*}
we observe
\beq\label{eq112}
\left|\sum_{i=1}^N(\bY_i-\bth)(\bY_i-\bth)'\right|=|\bS+N\bA|,
\eeq where $\bA=(\bth-\bar{\bY})(\bth-\bar{\bY})'$.

However, by taking advantage from Corollary A.3.1 of Anderson (2003) we reach the point that
\beq\label{eq113} |\bS+N\bA|&=&\left|\begin{array}{cc}
     S                               &   -\sqrt{N}(\bth-\bar{\bY})\\
     \sqrt{N}(\bth-\bar{\bY})'       &     1 \end{array}\right|\nonumber\\
&=&\left|\begin{array}{cc}
     1                               &   \sqrt{N}(\bth-\bar{\bY})'\\
     -\sqrt{N}(\bth-\bar{\bY})       &     S \end{array}\right|\nonumber\\
&=&|\bS|\{1+N(\bth-\bar{\bY})'\bS^{-1}(\bth-\bar{\bY})\} \eeq
Therefore, by making use of equations \eqref{eq111}-\eqref{eq113} we come to realize the following formula
\begin{eqnarray*}
f(\bth|\bY)&=&c(N,p)|\bS|^{-\frac{N}{2}}\{1+N(\bth-\bar{\bY})'\bS^{-1}(\bth-\bar{\bY})\}^{-\frac{N}{2}},
\end{eqnarray*}
where
\begin{eqnarray*}
c(N,p)&=&\left\{\int_{\bth\in\Re^p}|\bS|^{-\frac{N}{2}}\{1+N(\bth-\bar{\bY})'\bS^{-1}(\bth-\bar{\bY})\}^{-\frac{N}{2}}d\bth\right\}^{-1}\\
&=&|\bS|^{\frac{N}{2}}\left\{\frac{\left[\pi(N-p)\right]^{\frac{n}{2}}\Gamma\left(\frac{N-p}{2}\right)|\bS|^{\frac{1}{2}}}
{\Gamma\left(\frac{N}{2}\right)\left[N(N-p)\right]^{\frac{p}{2}}}
\right\}^{-1}\\
&=&\frac{N^{\frac{p}{2}}\Gamma\left(\frac{N}{2}\right)|\bS|^{\frac{N-1}{2}}}
{\pi^{\frac{n}{2}}(N-p)^{\frac{N-p}{2}}\Gamma\left(\frac{N-p}{2}\right)}.
\end{eqnarray*}
This would prove our claim.\hfill$\blacksquare$

Throughout this paper, we shall also assume that the loss function is
given by \beq\label{eq114}
L(\hat{\bth};\bth)=N(\hat{\bth}-\bth)'\bSi^{-1}(\hat{\bth}-\bth)
\eeq
for any estimator $\hat{\bth}$ of $\bth$.\\
It has been fully known that the Bayes estimator of $\bth$ with respect
to the loss \eqref{eq114} is the posterior mean (Proposition 2.5.1, Robert, 2001) given by
\begin{eqnarray}\label{eq115}
\hat{\bth}&=&\bar{\bY}.
\end{eqnarray}
As it can be realized from the estimator given by \eqref{eq115}, the Bayes estimator, reduces to the sample mean, under the setup presented above. So there is no need to deal with the Bayesian aspects of $\hat{\bth}$, and along the paper, we in fact concern sample mean rather than the Bayes estimator.

Then, from ECD properties (see Fang et al., 1990) we have
\beq\label{eq116} \hat{\bth}\sim \E_p(\bth,N^{-1}\bSi,g). \eeq
Under classical viewpoint, we devote a general class of Stein-type
shrinkage estimators to the estimator $\hat{\bth}$, given by
\beq\label{eq117}
\bd_r(\hat{\bth})=\left[1-\frac{r\left(\hat{\bth}^\prime\bS^{-1}\hat{\bth}\right)}
{\hat{\bth}^\prime\bS^{-1}\hat{\bth}}\right]\hat{\bth}, \eeq where
$r:[0,\infty)\rightarrow[0,\infty)$ is an absolutely continuous
function.

Furthermore, $r\in \mathcal{S}(\mathbb{R}_+,\mu)$, (the Schwartz space or space of rapidly decreasing functions on $\mathbb{R}_+$ with the measure $\mu$) where
\begin{eqnarray*}
\mathcal{S}(\mathbb{R}_+,\mu)=\left\{r\in\mathcal{C}^\infty(\mathbb{R}_+,\mu):\|r\|_{\alpha,\beta}<\infty\quad
\forall\;\alpha,\beta\right\},
\end{eqnarray*}
$\alpha$ and $\beta$ are indices, $\mathcal{C}^\infty(\mathbb{R}_+,\mu)$ is the set of all smooth functions from $\mathbb{R}_+$ to $\mathbf{C}$ (the set of all complex numbers) and
\begin{eqnarray*}
\|r\|_{\alpha,\beta}=\|x^\alpha\mathrm{D}^\beta r\|_\infty=\sup\{|x^\alpha\mathrm{D}^\beta r(x)|:x\in \mbox{domain}\;\mbox{of}\; r\}.
\end{eqnarray*}
Here $\mathrm{D}^\beta$ stands for $\beta^{th}$ derivative of $r$.
See Folland (1999) for more details.

The latter condition plays strategic position in gaining main result. Note that for every function such as $r(.)$ belongs to $\mathcal{S}(\mathbb{R}_+,\mu)$, we have
\beq
\int_0^\infty r'(x)d\mu(x)&<&\infty,\\
\int_0^\infty r^2(x)d\mu(x)&<&\infty,
\eeq
More interesting that the Schwartz space is dense in the space of all functions satisfy the above conditions in (1.18) and (1.19).

The objective of this study is to construct conditions on $r(.)$ under which
$\bd_r(\hat{\bth})$ performs better than $\hat{\bth}$ in the
sense of having smaller risk w.r.t. the loss
function given by \eqref{eq114}.

This study is highly motivated by the work of Srivastava and
Bilodeau (1989). They chewed over a similar class of estimators to
\eqref{eq117}, substituting the function $r(.)$ with a constant
under classical decision theory. Although, as noted above, considering Bayesian
point of view does not offer substantial generality, taking vague prior, over the work
of Srivastava and Bilodeau (1989), because of \eqref{eq116}, the class specified in
\eqref{eq117} contains the class which was previously stated as a special
case.

\setcounter{equation}{0}
\section{Risk Derivations}
In this section, we give some lemmas to evaluate the risk function
of $\bd_r(\hat{\bth})$. Provided that if all expectations exist,
we deserve the following Lemma.
\begin{lem}
If $\bx\sim\N_p(\bth,\alpha\bSi)$, $\alpha>0$,
$\bSi\in\mathcal{S}(p)$ is independent of $\bS\sim
W_p(\beta\bSi,n)$, $\beta>0$, $n=N-1$, then
\begin{eqnarray*}
E\left[\frac{\bx'\bSi^{-1}(\bx-\bth)r\left(\bx'\bS^{-1}\bx\right)}{\bx'\bS^{-1}\bx}\right]&=&
\beta\alpha(n-p+1)\bigg\{(p-2)E\left[\frac{r(\bx'\bS^{-1}\bx)}{\bx'\bSi^{-1}\bx}\right]\\
&&+2E\left[r'(\bx'\bS^{-1}\bx)\right]\bigg\}
\end{eqnarray*}
\vspace{-.3cm}and\vspace{-.3cm}
\begin{eqnarray*}
E\left[\frac{\bx'\bSi^{-1}\bx
r^2\left(\bx'\bS^{-1}\bx\right)}{\left(\bx'\bS^{-1}\bx\right)^2}\right]
&=&\beta^2(n-p+1)(n-p+3)E\left[\frac{r^2(\bx'\bS^{-1}\bx)}{\bx'\bSi^{-1}\bx}\right]
\end{eqnarray*}
\end{lem}
\noindent\textbf{Proof:} Let $\by=\bSi^{-\frac{1}{2}}\bx$,
$\bth^*=\bSi^{-\frac{1}{2}}\bth$, and
$\bA=\beta^{-1}\bSi^{-\frac{1}{2}}\bS\bSi^{-\frac{1}{2}}$, where
$\bSi^{-\frac{1}{2}}$ is a symmetric square root of $\bSi^{-1}$.
From independency of $\bx$ and $\bS$, $\by\sim
\N_p(\bth^*,\alpha\bI_p)$ is independent of both $\bA\sim
W_p(\bI_p,n)$ and
$\frac{\by'\by}{\by'\bA^{-1}\by}\sim\chi^2_{n-p+1}$. Also note
that $F=\bx'\bS^{-1}\bx=\beta^{-1}\by'\bA^{-1}\by$. Therefore
using Stein's (1981) identity we get
\begin{eqnarray*}
E\left[\frac{\bx'\bSi^{-1}(\bx-\bth)r\left(\bx'\bS^{-1}\bx\right)}{\bx'\bS^{-1}\bx}\right]&=&
E\left[\frac{\by'(\by-\bth^*)r(F)}{F}\right]\\
&=&
E\left[\frac{\by'(\by-\bth^*)r(F)}{\by'\by}\right]E\left[\frac{\by'\by}{F}\right]\\
&=&
\beta(n-p+1)E\left[\frac{\by'(\by-\bth^*)r(F)}{\by'\by}\right]\\
&=&\alpha\beta(n-p+1)\bigg\{(p-2)E\left[\frac{r(F)}{\by'\by}\right]+2E\left[r'(F)\right]
\bigg\}.
\end{eqnarray*}
Similarly
\begin{eqnarray*}
E\left[\frac{\bx'\bSi^{-1}\bx
r^2\left(\bx'\bS^{-1}\bx\right)}{\left(\bx'\bS^{-1}\bx\right)^2}\right]&=&
E\left[\frac{\by'\by r^2(F)}{F^2}\right]\\
&=&\beta^2E\left[\frac{r^2(F)}{\by'\by}\right]
E\left[\left(\frac{\by'\by}{\by'\bA^{-1}\by}\right)^2\right]\\
&=&\beta^2(n-p+1)(n-p+3)E\left[\frac{r^2(F)}{\by'\by}\right].
\end{eqnarray*}
\hfill$\blacksquare$
\begin{lem}
The risk function of the estimator $\bd_r(\hat{\bth})$ w.r.t. the loss function \eqref{eq112} is give by
\begin{eqnarray*}
R(\bd_r(\hat{\bth});\bth)&=&R(\hat{\bth};\bth)-4(N-p)\int_0^\infty
E\left[r'\left(\hat{\bth}^\prime
\bS^{-1}\hat{\bth}\right)\bigg|\;t\right]t^{-2}\W(dt)\\
&&+\int_0^\infty E\bigg\{\frac{(N-p)r\left(\hat{\bth}^\prime
\bS^{-1}\hat{\bth}\right)}{t\hat{\bth}^\prime
\left(t^{-1}\bSi\right)^{-1}\hat{\bth}}\\
&&\times\left[N(N-p+2)r\left(\hat{\bth}^\prime
\bS^{-1}\hat{\bth}\right)-2(p-2)\right]\bigg|\;t\bigg\}\W(dt)
\end{eqnarray*}
\end{lem}
\noindent\textbf{Proof:} As far as the representation
\eqref{eq14} is concerned, it is possible to continue this way
\beq\label{eq21}
R(\bd_r(\hat{\bth});\bth)&=&NE\left[(\bd_r(\hat{\bth})-\bth)'\bSi^{-1}(\bd_r(\hat{\bth})-\bth)\right]\nonumber\\
&=&R(\hat{\bth};\bth)-2N\int_0^\infty
E\left[\frac{\hat{\bth}^\prime\bSi^{-1}(\hat{\bth}-\bth)r\left(\hat{\bth}^\prime
\bS^{-1}\hat{\bth}\right)}{\hat{\bth}^\prime\bS^{-1}\hat{\bth}}\bigg|\;t\right]\W(dt)\nonumber\\
&&+N\int_0^\infty
E\left[\frac{\hat{\bth}^\prime\bSi^{-1}\hat{\bth}
r^2\left(\hat{\bth}^\prime
\bS^{-1}\hat{\bth}\right)}{\left(\hat{\bth}^\prime\bS^{-1}\hat{\bth}\right)^2}\bigg|\;t\right]\W(dt).
\eeq But from $\be_i|t\sim \N_p(\0,t^{-1}\bSi)$, it is concluded
that using \eqref{eq116}, $\hat{\bth}|t\sim
\N_p(\bth,t^{-1}N^{-1}\bSi)$ is independent of $\bS|t\sim
W_p(t^{-1}\bSi,n)$. Consequently, by making use of Lemma 2.1 for
$\alpha=(tN)^{-1}$ and $\beta=t^{-1}$ we get
\begin{eqnarray*}
E\left[\frac{\hat{\bth}^\prime\bSi^{-1}(\hat{\bth}-\bth)r\left(\hat{\bth}^\prime
\bS^{-1}\hat{\bth}\right)}{\hat{\bth}^\prime\bS^{-1}\hat{\bth}}\bigg|\;t\right]&=&\frac{N-p}{Nt^2}\bigg\{(p-2)
E\left[\frac{r\left(\hat{\bth}^\prime
\bS^{-1}\hat{\bth}\right)}{\hat{\bth}^\prime\bSi^{-1}\hat{\bth}}\bigg|\;t\right]\\
&&+2E\left[r'\left(\hat{\bth}^\prime
\bS^{-1}\hat{\bth}\right)\bigg|\;t\right]\bigg\},\\
E\left[\frac{\hat{\bth}^\prime\bSi^{-1}\hat{\bth}
r^2\left(\hat{\bth}^\prime
\bS^{-1}\hat{\bth}\right)}{\left(\hat{\bth}^\prime\bS^{-1}\hat{\bth}\right)^2}\bigg|\;t\right]&=&
\frac{(N-p)(N-p+2)}{t^2}
E\left[\frac{r^2\left(\hat{\bth}^\prime\bS^{-1}\hat{\bth}\right)}{\hat{\bth}^\prime\bSi^{-1}\hat{\bth}}\bigg|\;t\right].
\end{eqnarray*}
After all, substituting the above expressions in \eqref{eq21}, completes the proof.\hfill$\blacksquare$

\setcounter{equation}{0}
\section{Main Results}
In this section, we demonstrate the minimaxity of the estimator
$\bd_r(\hat{\bth})$, under some mild conditions made on the
function $r(.)$. Also we give conditions under which $\bd_r(\hat{\bth})$ dominates a James-Stein type shrinkage estimator.
\begin{thm}
Assume in the model \eqref{eq11}, $\be_i\sim \E_p(\0,\bSi,g)$.
Then w.r.t. the loss function \eqref{eq114}, the estimator
$\bd_r(\hat{\bth})$ is minimax in the sub-class \eqref{eq17}, providing
\begin{enumerate}
\item[\textbf{(i)}]$r$ is non-decreasing
\item[\textbf{(ii)}]$r\leq\frac{2(p-2)}{N(N-p+2)}$
\end{enumerate}
\end{thm}
\noindent\textbf{Proof:} The estimator $\hat{\bth}$ given by
\eqref{eq115} is minimax. Therefore, in order to show that
$\bd_r(\hat{\bth})$ is minimax it is enough to show that
$R(\hat{\bth};\bth)-R(\bd_r(\hat{\bth});\bth)\geq0$.
%Simplifying
%the terms the risk difference $\mathcal{D}$ is non-negative if
%\begin{eqnarray}\label{eq31}
%\mathcal{D}=2E\bigg\{F^{-1}r(F)\bx'\bSi^{-1}(\bx-\bth)-F^{-2}r^2(F)\bx\bSi^{-1}\bx\bigg\}\geq0,
%\end{eqnarray}
%where $F=\bx'\bS^{-1}\bx$. Let $\by=\bSi^{-\frac{1}{2}}\bx$,
%$\bth^*=\bSi^{-\frac{1}{2}}\bth$, and
%$\bA=\beta^{-1}\bSi^{-\frac{1}{2}}\bS\bSi^{-\frac{1}{2}}$, where
%$\bSi^{-\frac{1}{2}}$ is a symmetric square root of $\bSi^{-1}$.
%From independency of $\bx$ and $\bS$, $\by\sim
%\N_p(\bth^*,\alpha\bI_p)$ is independent of both $\bA\sim
%W_p(\bI_p,n)$ and
%$G=\frac{\by'\by}{\by'\bA^{-1}\by}\sim\chi^2_{n-p+1}$. Also note
%that $F=\beta^{-1}\by'\bA^{-1}\by$. Thus \eqref{eq31} changes to
%\begin{eqnarray}\label{eq32}
%\mathcal{D}=2E\bigg\{F^{-1}r(F)\by'(\by-\bth^*)-F^{-2}r^2(F)\by'\by\bigg\}\geq0.
%\end{eqnarray}
%It can be easily checked that from Stein's (1981) identity the
%condition \eqref{eq32} simplifies to
%\begin{eqnarray}
%\mathcal{D}=E\bigg\{2\alpha(p-2)F^{-1}r(F)-F^{-2}r^2(F)\|\by\|^2+4\alpha(\partial/\partial
%F)r(F)\bigg\}\geq0.
%\end{eqnarray}
%Thus $\mathcal{D}\geq0$ if
%\begin{eqnarray}
%E\bigg\{2\alpha(p-2)F^{-1}r(F)-F^{-2}r^2(F)\|\by\|^2+4\alpha(\partial/\partial
%F)r(F)\bigg\}\geq0
%\end{eqnarray}
%Therefore using Stein's (1981) identity we get
But from Lemma 2.2 we have
$R(\bd_r(\hat{\bth});\bth)-R(\hat{\bth};\bth)=\A+\B$, where
\begin{eqnarray*}
\A&=&-4(N-p)\int_0^\infty
E\left[r'\left(\hat{\bth}^\prime
\bS^{-1}\hat{\bth}\right)\bigg|\;t\right]t^{-2}\W(dt)\nonumber\\
\B&=&+\int_0^\infty E\bigg\{\frac{(N-p)r\left(\hat{\bth}^\prime
\bS^{-1}\hat{\bth}\right)}{t\hat{\bth}^\prime
\left(t^{-1}\bSi\right)^{-1}\hat{\bth}}\nonumber\\
&&\times\left[N(N-p+2)r\left(\hat{\bth}^\prime
\bS^{-1}\hat{\bth}\right)-2(p-2)\right]\bigg|\;t\bigg\}\W(dt)\nonumber\\
\end{eqnarray*}
Whereof $r(.)$ is non-decreasing, $r'\left(\hat{\bth}^\prime
\bS^{-1}\hat{\bth}\right)\geq0$. Also following Srivastava and
Bilodeau (1989), we have
\begin{eqnarray*}
&&\int_0^\infty
E\left[r'\left(\hat{\bth}^\prime
\bS^{-1}\hat{\bth}\right)\bigg|\;t\right]t^{-2}\W(dt)\\
&&=\int_0^\infty E\left[r'\left(\hat{\bth}^\prime
\bS^{-1}\hat{\bth}\right)\bigg|\;t\right]t^{-2}\W^+(dt)-
\int_0^\infty E\left[r'\left(\hat{\bth}^\prime
\bS^{-1}\hat{\bth}\right)\bigg|\;t\right]t^{-2}\W^-(dt)\\
&&\geq0.
\end{eqnarray*}
Therefore $\A\leq0$ and by making use of (1.18), $\A<\infty$.
Taking $r\leq\frac{2(p-2)}{N(N-p+2)}$, $\B\leq0$ is achieved for finite $\B$, which completes the proof.
But for demonstrating that $\B<\infty$, it is sufficient to show that
\beq\label{eq31}
&(i)&\int_0^\infty t^{-1}E\left(\frac{1}{N\hat{\bth}^\prime
\left(t^{-1}\bSi\right)^{-1}\hat{\bth}}\bigg|\;t\right)\W(dt)<\infty,\\
&(ii)&\int_0^\infty t^{-1}E\left\{\frac{r^2\left(\hat{\bth}^\prime
\bS^{-1}\hat{\bth}\right)}{N\hat{\bth}^\prime
\left(t^{-1}\bSi\right)^{-1}\hat{\bth}}\bigg|\;t\right\}\W(dt)<\infty.\nonumber
\eeq Note that for a fixed $t$, $N\hat{\bth}^\prime
\left(t^{-1}\bSi\right)^{-1}\hat{\bth}$ has non-central chi-square
distribution with $p$ d.f. and non-centrality parameter
$Nt\bth'\bSi^{-1}\bth$. In conclusion, the
\begin{eqnarray*}
E\left[\frac{1}{N\hat{\bth}^\prime \left(t^{-1}\bSi\right)^{-1}\hat{\bth}}\right]\leq
E\left[\frac{1}{\chi^2_p}\right]=\frac{1}{p-2}
\end{eqnarray*}
is observed and \eqref{eq31} (i) is followed by \eqref{eq17} (ii)-(iii).

On the other hand, using the covariance inequality (see Lemma 6.6 page 370 of Lehmann and Casella, 1998) and equation (1.19)
\begin{eqnarray*}
E\left\{\frac{r^2\left(\hat{\bth}^\prime
\bS^{-1}\hat{\bth}\right)}{N\hat{\bth}^\prime
\left(t^{-1}\bSi\right)^{-1}\hat{\bth}}\bigg|\;t\right\}\leq
E\left\{r^2\left(\hat{\bth}^\prime
\bS^{-1}\hat{\bth}\right)\bigg|\;t\right\}E\left\{\frac{1}{N\hat{\bth}^\prime
\left(t^{-1}\bSi\right)^{-1}\hat{\bth}}\bigg|\;t\right\}<\infty,
\end{eqnarray*}
Therefore \eqref{eq31} (ii) is satisfied by \eqref{eq17} (ii)-(iii).\hfill$\blacksquare$

In the following, we develop necessary conditions for the shrinkage estimator $\bd_r(\hat{\bth})$ to dominate the James-Stein type estimator given by
\begin{eqnarray}
\bd_{JS}(\hat{\bth})&=&\left[1-\frac{p-2}{\hat{\bth}^\prime\bS^{-1}\hat{\bth}}\right]\hat{\bth}.
\end{eqnarray}
The performance of this estimator is discussed in Srivastava and Bilodeau (1989) extensively. The way we derive the necessary conditions honesty is due to Maruyama and Strawderman (2005). But this approach has the following superiorities comparing to their analysis. (1) We study correlated errors with unknown covariance matrix while they considered uncorrelated case. (2) They derived the dominating result for multivariate normal, and we extend it for ECDs. Although the item (1) is completely different from that of uncorrelated, it is worthwhile to note that their conditions are robust under departures from normality assumptions. The following result is the same as Corollary 2.1. of Maruyama and Strawderman (2005). They could also find the class of admissible estimators under normal theory with identity covariance matrix.
\begin{thm}
Assume that the function $r(.)$ is bounded and absolutely continuous. Necessary conditions for an estimator $\bd_r(\hat{\bth})$ to dominate $\bd_{JS}(\hat{\bth})$ are that
\begin{enumerate}
\item[(i)] for every $\omega$, there exists $\omega_0$$(>\omega)$ such that $r'(\omega_0)\geq0$,
\item[(ii)] if $\omega r'(\omega)$ has a limiting value as $\omega$ approaches infinity, it must be $0$,
\item[(iii)] if $r(\omega)$ has a limiting value as $\omega$ approaches infinity and $\omega r'(\omega)$ converges to $0$ as $\omega$ approaches infinity, the limit value for $r(\omega)$ must be $\frac{p-2}{N(N-p+2)}$.
\end{enumerate}
\end{thm}
\noindent\textbf{Proof:} Proof of $(i)$ directly follows from the
proof of Corollary 2.1. of Maruyama and Strawderman (2005). Now
consider using Lemma 2.2, one can directly obtain \beq
\Delta &=& R(\bd_{JS}(\hat{\bth});\bth)-R(\bd_r(\hat{\bth});\bth)\nonumber\\
       &=& 4(N-p)\int_0^\infty
E\left[r'\left(\hat{\bth}^\prime
\bS^{-1}\hat{\bth}\right)\bigg|\;t\right]t^{-2}\W(dt)\nonumber\\
       &&+\int_0^\infty E\bigg\{\frac{(N-p)(p-2)}{t\hat{\bth}^\prime
\left(t^{-1}\bSi\right)^{-1}\hat{\bth}}\nonumber\\
&&\times\left[N(N-p+2)(p-2)-2(p-2)\right]\bigg|\;t\bigg\}\W(dt)\nonumber\\
       &&-\int_0^\infty E\bigg\{\frac{(N-p)r\left(\hat{\bth}^\prime
\bS^{-1}\hat{\bth}\right)}{t\hat{\bth}^\prime
\left(t^{-1}\bSi\right)^{-1}\hat{\bth}}\nonumber\\
&&\times\left[N(N-p+2)r\left(\hat{\bth}^\prime
\bS^{-1}\hat{\bth}\right)-2(p-2)\right]\bigg|\;t\bigg\}\W(dt)\nonumber\\
       &=&4(N-p)\int_0^\infty
E\left[r'\left(\hat{\bth}^\prime
\bS^{-1}\hat{\bth}\right)\bigg|\;t\right]t^{-2}\W(dt)\nonumber\\
       &&-(N-p)\int_0^\infty E\left[-\frac{\left\{r\left(\hat{\bth}^\prime
\bS^{-1}\hat{\bth}\right)-(p-2)\right\}^2}{\hat{\bth}^\prime
\bSi^{-1}\hat{\bth}}\bigg|t\right]t^{-2}\W(dt)\nonumber\\
&=&4(N-p)\int_0^\infty E\left[r'\left(\bz'\bB^{-1}\bz\right)\bigg|\;t\right]t^{-2}\W(dt)\nonumber\\
       &&-(N-p)\int_0^\infty E\left[-\frac{\left\{r\left(\bz'\bB^{-1}\bz\right)-(p-2)\right\}^2}{\bz'\bB^{-1}\bz}
       \frac{\bz'\bB^{-1}\bz}{\bz'\bz}\bigg|t\right]t^{-2}\W(dt)\nonumber\\
       &=&(N-p)\int_0^\infty E\left[\mathcal{G}_{r}(\bz'\bB^{-1}\bz)\bigg| t\right]t^{-2}\W(dt), \eeq where
$\bz=\bSi^{-\frac{1}{2}}\bar{\bY}$,
$\bB=\bSi^{-\frac{1}{2}}\bS\bSi^{-\frac{1}{2}}$ and \beq
\mathcal{G}_{r}(\omega)=-\frac{\left[r(\omega)-(p-2)\right]^2}{(n-p-1)\omega}+4r'(\omega).
\eeq For the proofs of $(ii)$ and $(iii)$, using the proofs of
$(ii)$ and $(iii)$ of Corollary 2.1. of Maruyama and Strawderman
(2005), it is enough to show that if $\bd_r(\hat{\bth})$ dominates
$\bd_{JS}(\hat{\bth})$, then, for every $\omega$, there exists
$\omega_0$$(>\omega)$ such that $\mathcal{G}_{r}(\omega_0)\geq0$.
In this case we follow the proof of Theorem 2.1. of Maruyama and
Strawderman (2005).

Suppose to the contrary that there exists $\omega_0$ such that $\mathcal{G}_{r}(\omega)<0$ for any $\omega\geq\omega_0$. Under the boundedness of $\mathcal{G}_{r}(.)$, there exists an $M$$(>0)$ such that $\mathcal{G}_{r}(\omega)\leq M$ for any $\omega$. Under the assumption of absolute continuity of $\mathcal{G}_{r}(.)$ there exists two points $(\omega_0<)$$\omega_1<\omega_2$ and $\epsilon$$(>0)$ such that $\mathcal{G}_{r}(\omega)<-\epsilon$ on $\omega\in[\omega_1,\omega_2]$. Using $M$ and $\epsilon$, we define $\mathcal{G}_{r,\epsilon}(\omega)$ as
\beq
\mathcal{G}_{r,\epsilon}(\omega)=\left\{\begin{array}{cc}
M         & \omega\leq\omega_0\\
0         & \omega_0<\omega<\omega_1\\
-\epsilon & \omega_1\leq\omega\leq\omega_2\\
0         & \omega>\omega_2
                                       \end{array}\right.
\eeq
The inequality $\mathcal{G}_{r,\epsilon}(\omega)\geq\mathcal{G}_{r}(\omega)$ for any $\omega$ and using equation (3.3) imply
\beq
\Delta&=&(N-p)\int_0^\infty E\left[\mathcal{G}_{r}(\bz'\bB^{-1}\bz)\bigg| t\right]t^{-2}\W(dt)\nonumber\\
&\leq& M(N-p)\int_0^\infty P_{\bth}\left(W\leq\omega_0\bigg|t\right)t^{-2}\W(dt)\nonumber\\
&&-\epsilon(N-p)\int_0^\infty P_{\bth}\left(\omega_1\leq W\leq\omega_2\bigg|t\right)t^{-2}\W(dt),
\eeq
where $W=\|\bX\|^2$ for $\bX=(t^{-1}\bSi)^{-\frac{1}{2}}\hat{\bth}$.\\
Based on the properties of the model under study, it can be realized that
\begin{eqnarray*}
&&\int_0^\infty P_{\bth}\left(W\leq\omega_0\bigg|t\right)t^{-2}\W(dt)\\
&&=\int_0^\infty
P_{\bth}\left(W\leq\omega_0\bigg|t\right)t^{-2}\W^+(dt)-
\int_0^\infty
P_{\bth}\left(W\leq\omega_0\bigg|t\right)t^{-2}\W^-(dt)\\
&&\geq0.
\end{eqnarray*}
This phenomenon is also valid for $\int_0^\infty P_{\bth}\left(\omega_1\leq W\leq\omega_2\bigg|t\right)t^{-2}\W(dt)$.

Now let $\ba$ be a fixed $p$-dimensional unit vector (see Fig. 1). Then the half plane $\{\bx:\ba'\bx\leq\sqrt{\omega_0}\}$ includes the $p$-dimensional hyper-ellipsoid $\{\bx:\|\bx\|^2\leq\omega_0\}$. For $\bth=(\sqrt{\omega_0}+\lambda)(t^{-1}\bSi)^{-\frac{1}{2}}\ba$, we have
\beq
P_{\bth}\left(W\leq\omega_0\bigg|t\right)&<&\int_{\ba'\bx\leq\sqrt{\omega_0}}\frac{|t^{-1}\bSi|^{-\frac{1}{2}}}
{(2\pi)^{\frac{p}{2}}}\;\exp\left(-\frac{\|\bx-\bth\|^2}{2}\right)d\bx\nonumber\\
&\leq&\exp(\lambda\sqrt{\omega_0})\exp\left(-\frac{\|\bth\|^2}{2}\right)\nonumber\\
&&\times\int_{\ba'\bx\leq\sqrt{\omega_0}}\frac{|t^{-1}\bSi|^{-\frac{1}{2}}}
{(2\pi)^{\frac{p}{2}}}\;\exp\left(-\frac{\|\bx\|^2}{2}+\sqrt{\omega_0}\ba'\bx\right)d\bx\nonumber\\
&\leq&\exp(\lambda\sqrt{\omega_0})\exp\left(-\frac{\|\bth\|^2}{2}+\frac{\omega_0}{2}\right).
\eeq
For $N=\{\bx:\omega_1\leq\|\bx\|^2\leq\omega_2,\sqrt{\omega_1}\leq\ba'\bx\leq\sqrt{\omega_2}\}$ and $\bth=(\sqrt{\omega_0}+\lambda)(t^{-1}\bSi)^{-\frac{1}{2}}\ba$, we get
\beq
P_{\bth}\left(\omega_1\leq W\leq\omega_2\bigg|t\right)&>&\int_N\frac{|t^{-1}\bSi|^{-\frac{1}{2}}}
{(2\pi)^{\frac{p}{2}}}\;\exp\left(-\frac{\|\bx-\bth\|^2}{2}\right)d\bx\nonumber\\
&\geq&\exp(\lambda\sqrt{\omega_1})\exp\left(-\frac{\|\bth\|^2}{2}\right)\nonumber\\
&&\times\int_N\frac{|t^{-1}\bSi|^{-\frac{1}{2}}}
{(2\pi)^{\frac{p}{2}}}\;\exp\left(-\frac{\|\bx\|^2}{2}+\sqrt{\omega_0}\ba'\bx\right)d\bx.
\eeq
By making use of the equations (3.6)-(3.8), we can obtain
\begin{eqnarray*}
\Delta\leq (N-p)\int_0^\infty c_1\exp\left(\sqrt{\omega_0}\lambda-\frac{\|\bth\|^2}{2}\right)
\bigg(1-c_2\exp\bigg[(\sqrt{\omega_1}-\sqrt{\omega_0})\lambda\bigg]\bigg)t^{-2}\W(dt),
\end{eqnarray*}
where $c_1=M\exp\left(\frac{\omega_0}{2}\right)$ and
\begin{eqnarray*}
c_2=\frac{\epsilon}{M}\;\exp\left(\frac{\omega_0}{2}\right)\int_N\frac{|t^{-1}\bSi|^{-\frac{1}{2}}}
{(2\pi)^{\frac{p}{2}}}\;\exp\left(-\frac{\|\bx\|^2}{2}+\sqrt{\omega_0}\ba'\bx\right)d\bx.
\end{eqnarray*}
Since $c_1$ and $c_2$ do not depend on $\lambda$, $\Delta$ is negative for sufficiently large $\lambda$. This completes the proof.\hfill$\blacksquare$\\
\input{epsf}
\begin{figure}
\epsfxsize=3in
\epsfysize=3in
\epsffile{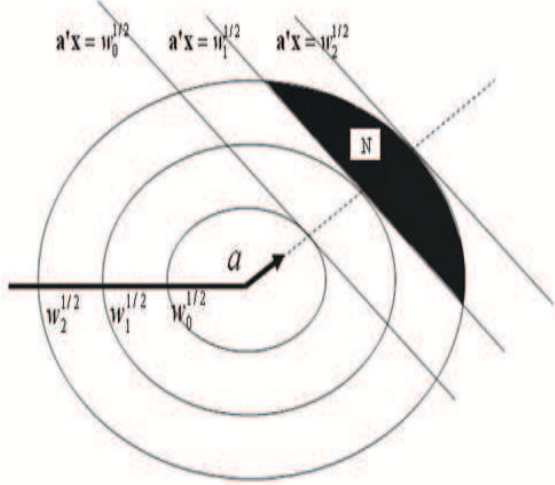}
\caption{\textbf{Graph of half planes}}
\end{figure}
Subsequently, we continue on giving an example of the function $r(.)$.

Let
\begin{eqnarray}\label{eq39}
r_*(x)=\frac{(p-2)b}{1+cx^{-1}},\hspace{.5cm}b=\frac{1}{N(N-p+2)}\quad\mbox{and}\quad c\in\mathbb{R}_+.
\end{eqnarray}
Then we have
\begin{eqnarray*}
0<r_*(x)&\leq&\frac{2(p-2)}{N(N-p+2)},\hspace{1cm}\mathrm{D}r_*(x)=\frac{b(p-2)^2}{(x+c)^2}>0,\\
\lim_{x\rightarrow\infty}r_*(x)&=&b(p-2),\hspace{1cm}\lim_{x\rightarrow\infty}x\mathrm{D}r_*(x)=0,
\end{eqnarray*}
which satisfy the conditions of Theorems 3.1 and 3.2

The resulting shrinkage estimator using the function $r_*(.)$ in \eqref{eq39}, is the generalized type of Alam and Thompson (1969) estimator given by
\begin{eqnarray*}
\bd_{r_*}(\hat{\bth})=\left\{1-\frac{r_*\left(\hat{\bth}^\prime\bS^{-1}\hat{\bth}\right)}
{\hat{\bth}^\prime\bS^{-1}\hat{\bth}}\right\}\hat{\bth}=\left\{1-\frac{(p-2)b}
{\hat{\bth}^\prime\bS^{-1}\hat{\bth}+c}\right\}\hat{\bth}.
\end{eqnarray*}
Also note that based on (1.18) and (1.19), the required conditions
of the Schwartz space, for this example, are
$b(p-2)^2E(X+c)^{-2}<\infty$ and
$b^2(p-2)^2E\left(\frac{X}{X+c}\right)^{2}<\infty$, which
summarizes to the sole condition
$E\left(\frac{X}{X+c}\right)^{2}<\infty$.

\setcounter{equation}{0}
\section{Conclusions}
In this paper, we utilized a broad class of Stein-type estimators which outperformed the consistent estimator of the mean of an elliptically contoured model.
 It is worthwhile to note that the minimaxity conditions are identical to that obtain under normal assumptions. Hence, those are robust with respect to departures from normality. Moreover, Bayesian perpective dose not offer systematic generality over classical approaches taking flat prior information. The class of estimators considered in Srivastava and Bilodeau (1989) is broaden into a more general shrinkage estimators; and as a result, this work dominates series of Brandwein's and Berger's papers. To the best of my knowledge, it is not simple to prove the admissibility of the class of Bayes shrinkage estimator $\bd_r(\hat{\bth})$ under elliptical symmetry and there exists no study in ECDs when the covariance matrix in unknown. But one may demonstrate it through taking the harmonic prior $\|\bth\|^{2-p}$ for $\pi(\bth)$ in (1.9) which leaves for further research. In this case, one may follow the work of Maruyama (2004) under the integral representation of elliptical models in (1.4). In this case, the work of Fourdrinier et al. (2003) has some interesting features.

%\section{Acknowledgements}
%The author would like to thank an anonymous referee for many
%valuable and constructive comments that led to vast improvements
%of the paper.

\section*{References}

\baselineskip=12pt
\def\ref{\noindent\hangindent 25pt}

\ref Alam, K. and Thompson, J. R., (1969). Locally averaged risk, {\em Ann. Inst. Statist. Math.}, {\bf21}, 457-469.

\ref Anderson, T. W., (2003). {\em An introduction to multivariate
statistical analysis}, 3rd ed., John Wiley and Sons, New York.

\ref Anderson, T. W. and Fang, K. T., (1990). Inference in multivariate elliptically contoured distribution based on maximum likelihood in "{\em Statistical Inference in Elliptically Contoured and Related Distribution}", (K. T. Fang and T. W. Anderson Ed.), 201-216. Allerton Press, New York.

\ref Anderson, T. W., Fang, K. T. and Hsu, H., (1986).
Maximum-likelihood estimates and likelihood-ratio criteria for
multivariate elliptically contoured distributions, {\em The
Canadian J. Statist.}, {\bf 14}, 55--59.

\ref Berger, J., (1975). Minimax estimation of location vectors for a wide class of densities, {\em Ann. Statist.}, {\bf3}, 1318-1328.

\ref Berger, J., (1976). Minimax estimation of a multivariate normal mean under arbitrary quadratic loss, {\em J. Mult. Anal.}, {\bf6}, 256-264.

\ref Box, G. E. P. and Tiao, G. C. (1992) {\em Bayesian Inference in Statistical Analysis}, John Wiley, New York.

\ref Brandwein, A. C., (1979). Minimax estimation for the mean of spherically symmetric distributions under general quadratic loss, {\em J. Mult. Anal.}, {\bf9}, 579-588.

\ref Brandwein, A. C. and Strawderman, W. E., (1978). Minimax estimation of location parameters of spherically symmetric unimodal distributions under quadratic loss, {\em Ann. Statist.}, {\bf6}, 377-416.

\ref Cambanis, S., Huang, S. and Simons, G., (1981). On the theory of elliptically contured distributions, {\em J. Mult. Anal.}, {\bf11}, 368–385.

\ref Cellier, D., Fourdrinier, D. and Robert, C., (1989). Robust shrinkage estimators of the location parameter for elliptically symmetric distributions, {\em J. Mult. Anal.}, {\bf29}, 39-52.

\ref Chu, K. C., (1973), Estimation and decision for linear systems with
elliptically random process. {\em IEEE Trans. Autom. Cont.}, 18,
499-505.

\ref Das Gupta, S., Eaton, M. L., Olkin, I., Perlman, M., Savage, L. J. and Sobel, M., (1972). Inequalities on the probability content of convex regions for elliptically contoured distributions, {\em Proc. Sixth Berkeley Symp. Math. Statist. Prob. vol. II. Prob. Theo.}, University of California Press, Berkeley, California, 241–265.

\ref Debnath, L. and Bhatta, D., (2007). {\em Integral Transforms and Their Applications}, Chapman and Hall, London, New York.

\ref Fang, K. T., Kotz, S. and Ng, K. W., (1990). {\em Symmetric
Multivariate and Related Distributions}, Chapman and Hall, London,
New York.

\ref Fang, K. T. and Zhang, Y., (1990), {\em Generalized Multivariate Analysis}, Springer, Beijing.

\ref Folland, G. B., (1999). {\em Real Analysis: Modern Techniques and Their Applications}, 2nd ed., John Wiley and Sons, New York.

\ref Fourdrinier, D., Strawderman, W. E. and Wells, M. T., (2003). Robust shrinkage estimation for elliptically
symmetric distributions with unknown covariance matrix, {\em J. Mult. Anal.}, {\bf85}, 24-39.

\ref Gupta, A. K. and Varga, T., (1993). {\em Elliptically
Contoured Models in Statistics}, Kluwer Academic Press.

\ref Jeffreys, H., (1961). {\em Theory of Probability}, Oxford: Clarendon.

\ref Kibria, B. M. G. and Haq, M. S., (1999). Predictive inference for the elliptical linear model, {\em J. Mult. Anal}, {\bf68}, 235-249.

\ref Lehmann, E. L. and Casella, G. (1998). {\em Theory of Point
Estimation}, 2nd ed., Springer, New York.

\ref Maruyama, Y. (2004). Stein's idea and minimax admissible estimation of a multivariate normal mean, {\em J. Mult. Anal.}, {\bf88}, 320-334.

\ref Maruyama, Y. and Strawderman, W. E. (2005). Necessary conditions for dominating the James-Stein estimator, {\em Ann. Inst. Statist. Math.}, {\bf57}, 157-165.

\ref Muirhead, R. J., (1982). {em Aspect of Multivariate
Statistical Theory}, John Wiley, New York.

\ref Ng, V. M., (2002). Robust Bayesian Inference for Seemingly Unrelated Regressions with Elliptical Errors , {\em J. Mult. Anal.}, {\bf82}, 409-414.

\ref Robert, C. P., (2001). {\em The Bayesian Choice: from Decision-Theoretic Motivations to Computational Implementation}, 2nd ed., Springer, New York.

%\ref Srivastava, M. S. and Khatri, C. G., (1979). {\em An Introduction to Multivariate Analysis}, North-Holland, Amsterdam.

\ref Srivastava, M. and Bilodeau, M., (1989). Stein estimation under elliptical distribution, {\em J. Mult. Annal.}, {\bf28}, 247-259.

\ref Zellner, A. (1971). {\em An Introduction to Bayesian Inference in Econometrics}, John Wiley, New York.

\end{document}